\documentclass[12pt]{article}

\usepackage{amsmath}
\usepackage{amsfonts}
\usepackage{pgfplots}
\usepackage{amssymb}
\usepackage[utf8]{inputenc}
\usepackage{tikz}
\usetikzlibrary{positioning}
\usetikzlibrary{matrix,arrows,decorations.pathmorphing}
\usetikzlibrary{patterns}
\usetikzlibrary{through,calc,3d}
\usetikzlibrary{plotmarks}
\usepackage{varioref}
\usepackage{subfigure}                  

\newtheorem{theorem}{\bf Theorem}[section]

\newcommand{\pfrac}[2]{\ensuremath{\dfrac{\partial #1}{\partial #2}}}
\newcommand{\inner}[2]{\big< #1 ,  #2 \big>}
\newcommand{\bm}[1]{\ensuremath{\mathbf{#1}}}
\newcommand{\bs}[1]{\ensuremath{\boldsymbol{#1}}}
\newcommand{\ol}[1]{\ensuremath{\overline{#1}}}
\newcommand{\bxi}{\boldsymbol{\xi}}
\newcommand{\Wx}{\mathbf{W}_{\mathbf{x}}}
\newcommand{\Wxi}{\mathbf{W}_{\boldsymbol{\xi}}}
\newcommand{\lp}{\left(}
\newcommand{\rp}{\right)}

\begin{document}

\title{A Scalable Observation-Driven Time-Dependent Basis  for a Reduced Description of  Transient  Systems}

\author{Hessam Babaee$^{1}$
\thanks{{h.babaee@pitt.edu},
 }\\
 $^{1}$Department of Mechanical Engineering and Materials Science\\
University of Pittsburgh, 3700 O'Hara St., Pittsburgh, PA 15261, USA
}

\date{}
\maketitle
\begin{abstract}
We present a variational principle for the extraction of a time-dependent orthonormal basis from random realizations of transient systems. The optimality condition of the variational principle leads to a closed-form evolution equation for the orthonormal basis and its coefficients. The extracted modes are associated with the most transient subspace of the system, and they provide a reduced description of the transient dynamics that may be used for reduced-order modeling, filtering and prediction. The presented method is matrix-free and relies only on the observables of the system and ignores any information about the underlying system. In that sense, the presented reduction is purely observation-driven and may be applied to systems whose models are not known. The presented method has linear computational complexity and memory storage requirement with respect to the number of observables and the number of random realizations. Therefore, it may be used for a large number of observations and samples. The effectiveness of the proposed method is tested on three examples: (i) stochastic advection equation, (ii) a reduced description of transient instability of Kuramoto-Sivashinsky, and (iii) a transient vertical jet governed by incompressible Navier-Stokes equation. In these examples, we contrast the performance of the time-dependent basis versus static basis such as proper orthogonal decomposition, dynamic mode decomposition and polynomial chaos expansion.
\end{abstract}



 \section{Introduction}
The ability to discover intrinsic structure from data has been called the
\emph{the fourth paradigm of scientific discovery} \cite{HTT09}.
Dimension reduction  enables such discoveries and clearly belongs to the fourth paradigm. 
The goal of dimension reduction is to map high-dimensional data into a meaningful 
representation of reduced dimensionality.  
As a result, dimension reduction has been widely used in classification, visualization, and 
compression of high-dimensional data \cite{MPH07}. 
In the scientific and engineering applications dimension reduction has been used 
to build reduced order models to propagate uncertainty, control  and optimization \cite{Noack:2016aa,marzouk2009stochastic}. 

This paper focuses on  dimension reduction of vector-valued observables  $\bm{x}(t;\omega) \in \mathbb{R}^n$ generated by time-dependent stochastic system where $\omega$ denotes a random observation.  
A linear reduction model seeks to approximate the target $\bm{x}$ with
\begin{equation}\label{Eq:Lin_Rduc}
\bm{x} = \bm{U} \bm{z} + \epsilon
\end{equation}
where $\bm{U} \in \mathbb{R}^{n\times r}$ is a deterministic matrix whose columns can be regarded as \emph{modes} and $\bm{z} \in \mathbb{R}^r$ are the reduced stochastic variable, where $r$ is the reduction size and $\epsilon$ is the reduction error. Both $\bm{U}$ and $\bm{z}$ are in general time-dependent.  In the view of 
Equation \ref{Eq:Lin_Rduc}, polynomial chaos expansion (PCE) can be seen as a special case where $\bm{z}$ represented by polynomial chaos – a \emph{static} (i.e. time-independent) basis – while $\bm{U}=\bm{U}(t)$ is time-dependent. PCE has reached widespread popularity especially for steady-state and 
parabolic stochastic partial differential equations (SPDE) \cite{XK02}.  However, for time-dependent dynamical
 systems with transient instabilities  the number of PCE terms required for a given level of accuracy must rapidly increase \cite{BM12}. 
 Similarly, it was shown in \cite{WK06} that the reduction size of  PCE  must increase 
 in time to maintain a given accuracy. 
On the other hand,  data-driven methods 
 such as proper  orthogonal decomposition (POD) and dynamic mode decomposition (DMD) belong to a class of methods
 where the modes $\bm{U}$ are static and they are computed from numerical/experimental observations of the dynamical system and $\bm{z}(t; \omega)$
 is obtained by often projection of the data to these modes \cite{S10}.  For highly transient systems, DMD and POD require large number of modes to capture transient time behavior \cite{KBBJ16,CM04}. This is partly because  both DMD and POD decompositions seek to obtain an optimal reduction of the data in a \emph{time-average} sense. In the case of DMD a best-fit linear dynamical system  matrix $\bm{A}$ is sought so that $\|\bm{x}^{k+1}-\bm{A}\bm{x}^k \|_2 $ is minimized across all snapshots, where $k$ shows the time snapshots. Similarly, the POD is an optimal basis for all snapshots.  However, transient instabilities, for example, by definition have a short life time. The effect of these transient instabilities   may be lost in a reduction based on time average.


More recently dynamically orthogonal (DO) \cite{SL09}  and bi-orthogonal (BO) \cite{CHZI13} expansions  have proven promising
in tackling time-dependent SPDEs. In  these reduction  techniques  both  $\bm{U}(t)$ and $\bm{z}(t; \omega)$ are \emph{dynamic} 
and they evolve with time. In DO/BO framework,  exact evolution equations for the evolution of $\bm{U}(t)$ and 
$\bm{z}(t; \omega)$ are derived.  In the case of SPDEs, the evolution equation for $\bm{U}(t)$ is a system of 
deterministic PDEs that are solved fully coupled  with  a system of stochastic ordinary differential equations
(SODE) and  an evolution equation for the mean equation. Although DO and BO algorithms lead to different evolution equations for the stochastic 
coefficient  and the spatial basis, it is shown in \cite{CSK14,Babaee:2017aa} that both these methods are equivalent, i.e. the BO and DO 
reductions span the same subspaces. In deterministic applications,   time-dependent basis have been used for dimension reduction of 
 deterministic Schrodinger equations by the Multi Configuration Time Dependent Hartree (MCTDH) method \cite{Beck:2000aa,Bardos:2003ab}, and later by Koch and Lubich \cite{Koch:2007aa}. Recently, optimally time-dependent reduction (OTD) was introduced \cite{Babaee_PRSA}, in which an evolution equation for a set of  orthonormal modes is found that capture the most instantaneously transient subspace of a dynamical system \cite{BFHS17}.


The motivation for this paper is to obtain a reduced description of transient dynamics that relies only on random observations --- i.e. an equation-free reduction. This is in contrast to model-dependent reduction techniques, for example DO and BO methodologies, that require solving high-dimensional ODE/PDE for the determination of the time-dependent modes. The model-dependence of these techniques prevents the application of these methods to complex systems whose underlying models are not known or suffer from physics deficiency due to multi-scale physics effects,  for example,  atmospheric sciences and weather forecasting,  computational neuroscience and complex biological systems.   Moreover,  the evolution equation for the DO/BO modes is tied to the  equation that governs the dynamical system at hand, which would require a case by case implementation.

The objective of this paper is to find a reduced description of transient dynamics that is observation driven and discards any information about the underlying model.  To this end,  we present a variational principle for the determination of the set of time-dependent orthonormal modes that capture the low-dimensional structure from observable of the stochastic system. The structure of the paper is as follows: In Section 2, we present the variational principle for the extraction of the time-dependent modes and their coefficients. In Section 3, we present three demonstration cases: (a) the stochastic advection equation, (b) Kuramoto-Sivashinsky equation, and (c) a transient vertical jet in quiescent background flow governed by incompressible Navier-Stokes equations. In Section 4, we present the conclusions of our work.

%


 \section{Methodology}
\subsection{Setup}
We consider  a stochastic dynamical system  in the form of:
\begin{equation}\label{eq:state}
\dot{\mathbf{x}} = \mathbf{f}(\mathbf{x},t;\boldsymbol \xi),
\end{equation}
where  $\mathbf{x}$ is the state space and the above dynamical system is finite-dimensional, i.e. $\mathbf{x} \in \mathbb{R}^m$, or  infinite dimensional, expressed, for example, by SPDEs and $\boldsymbol \xi =(\xi_1, \xi_2, \dots, \xi_d) \in \mathbb{R}^d$ is a $d$-dimensional
random space with the joint probability given by $\rho(\boldsymbol \xi)$.  
The dynamical system  may not be known and the proposed reduction relies only on  vector-valued observations of the above system. 
We base our formulation on  continuous-time  observations, represented by matrix   $\bm{T}(t) \in \mathbb{R}^{n\times s}$, whose  columns contain vector-valued observations of  the stochastic system at time $t$:
\begin{align*}
\bm{T}(t) &= \bigg[\bm{x}_1(t)\  \Big| \ \bm{x}_2(t) \ \Big| \  \dots \ \Big| \ \bm{x}_s(t)   \bigg]_{n\times s}.
\end{align*}
Each column of this matrix is obtained by first drawing a random sample of $\boldsymbol \xi$ with joint probability distribution of $\rho(\boldsymbol \xi)$ and solving the forward model given by   equation \ref{eq:state}.
Before we proceed we introduce the sample mean operator: 
\begin{equation*}
\mathbb{E}[\bm{T}] = \bm{T} \mathbf{w}_{\boldsymbol \xi},
\end{equation*}
where $\mathbf{w}_{\boldsymbol \xi} = [\mathbf{w}_{{\boldsymbol \xi}_1}, \mathbf{w}_{{\boldsymbol \xi}_2}, \dots, \mathbf{w}_{{\boldsymbol \xi}_s}]^T \in \mathbb{R}^{s\times 1}$ represents the sample weights. For example, for standard Monte Carlo samples, $\mathbf{w}_{\boldsymbol \xi}=\bm{1}/s$, where $\bm{1}\in \mathbb{R}^{s\times 1}$ is a column vector with all elements equal to 1. For cases where controlled observations of system \ref{eq:state} is possible, the probabilistic collocation strategy may be used for  sampling. In that case 
$\mathbf{w}_{\boldsymbol \xi}$ represents the associated probabilistic collocation weights. Similarly, $\bm{w}_{\bm{x}}=[\bm{w}_{{\bm{x}}_1},\bm{w}_{{\bm{x}}_2}, \dots, \bm{w}_{{\bm{x}}_n}]^T   \in \mathbb{R}^{n\times 1}$ 
represents the associated weights in the state space. For applications governed by SPDEs, $\mathbf{w}_{\mathbf{x}}$ may be chosen as space integration weights (e.g. quadrature weights), such that $\bm{w}_{\bm{x}}^T \bm{u}$ approximates $\int_{\Omega} u d\Omega$, where $\Omega$ is the spatial domain and $\bm{u}$ contains vector-valued observations of spatial function $u$. 
We assume that $\bm{T}(t)$ has zero mean, i.e. $\mathbb{E}[\bm{T}(t)]=0$. This is  accomplished by subtracting the ensemble mean  from all the realizations.  Let's define the inner product in both state and random spaces:
 \begin{align*}
 \inner{\bm{u}_i}{\bm{u}_j} &=   \bm{u}^T_i \Wx \bm{u}_j,\\
 \mathbb{E}[\bm{y}^T_i \bm{y}_j] &=   \bm{y}^T_i \Wxi \bm{y}_j,
 \end{align*}
 where $\Wx = \mbox{diag} (\mathbf{w}_{\mathbf{x}})$ and $\Wxi = \mbox{diag} (\mathbf{w}_{\boldsymbol \xi})$
 are matrices of size $n \times n$ and $s \times s$, respectively. Therefore, $\inner{\bm{u}_i}{\bm{u}_j}$ is an inner product of  $\bm{u}_i$ and $\bm{u}_j$ in the state space, and  $ \mathbb{E}[\bm{y}^T_i \bm{y}_j]$ is an inner product of $\bm{y}_i$ and $\bm{y}_j$ in the random space. Our goal is to build a low-rank approximation
 that can capture transient behavior in the stochastic data $\bm{T}(t)$. To this end, we seek a reduction in the form of:   
\begin{equation}
\bm{T}(t) = \bm{U}(t) \bm{Y}(t)^T + \mathcal{E}(t),
\end{equation}
with 
\begin{align*}
\bm{Y}(t) &= \bigg[\bm{y}_1(t)\  \Big| \ \bm{y}_2(t) \ \Big| \  \dots \ \Big| \ \bm{y}_r(t)   \bigg]_{s\times r}, \\
\bm{U}(t) &= \bigg[\bm{u}_1(t)\  \Big| \ \bm{u}_2(t) \ \Big| \  \dots \ \Big| \ \bm{u}_r(t)   \bigg]_{n\times r}, 
\end{align*}
where $\bm{Y}(t) \in  \mathbb{R}^{s\times r}$  and $\bm{Y}(t) = \big[\bm{y}_1(t)\  | \ \bm{y}_2(t) \ | \  \dots \ | \ \bm{y}_r(t)   \big]$ and $\bm{U}(t) \in  \mathbb{R}^{n\times r}$ is a set of orthonormal modes, i.e. $\inner{\bm{u}_i(t)}{\bm{u}_j(t)} = \delta_{ij}$ and $\mathcal{E}$ is the reduction error. In the proposed form $\bm{U}(t)$ represents a time-dependent low rank subspace which we will refer to as \emph{dynamic basis} and $\bm{Y}(t)$ are the stochastic coefficients of the dynamic basis and  we will refer to $\bm{Y}(t)$ as the \emph{stochastic coefficients}. In the following we present a variational principle that results in  closed-form evolution equations for the dynamic basis $\bm{U}(t)$ and the stochastic coefficients $\bm{Y}(t)$.
\subsection{Variational Principle}
In this section we present a variational principle from which the evolution equations of the low-rank approximation are derived. Our
goal is to build a low-rank approximation that can capture transient data. To this end, we first introduce  a weighted Frobenius norm: 
\begin{equation}\label{eq:Wgh_Forb}
\left\Vert \mathbb{E}\big[ \bm{T}  \big]^{2}  \right\Vert_F^{1/2}  = \bigg( \sum_{i=1}^n \sum_{j=1}^s \bm{w}_{{\bm{x}}_i} \bm{w}_{{\boldsymbol \xi}_j} \bm{T}_{i,j}^2 \bigg)^{1/2},
\end{equation}
where $\bm{T}_{i,j}$ is the entry of matrix $\bm{T}$ at the $i$th row and the $j$th column. 
Now we consider the variational principle in the form of:
\begin{equation}\label{eq:functional}
\mathcal{F}(\dot{\bm{U}}, \dot{\bm{Y}}) 
 =\left\Vert \mathbb{E}\bigg[ \frac{d (\bm{U} \bm{Y}^T)}{dt}  - \dot{\bm{T} } \bigg]^{2}  \right\Vert_F,
\end{equation}
 subject to the orthonormality condition of the modes, i.e. $\inner{\bm{u}_i(t)}{\bm{u}_j(t)} = \delta_{ij}$. The variational principle \ref{eq:functional} seeks to minimize the difference between the rate of change of all observations and the dynamic basis reduction, and the control parameters are $\dot{\bm{U}}$ and $\dot{\bm{Y}}$.  To incorporate the orthonormality condition , we observe that:
\begin{equation}
\inner{\dot{\bm{u}}_i}{\bm{u}_j} + \inner{\bm{u}_i}{\dot{\bm{u}}_j} = 0,
\end{equation}
and we denote $\bs{\phi}_{ij}(t): = \inner{\bm{u}_i}{\dot{\bm{u}}_j}$. It is clear from the above equation that $\bs{\phi}_{ij}(t)$ is a skew-symmetric matrix, i.e. $\bs{\phi}_{ij}(t) = -\bs{\phi}_{ji}(t)$.
We incorporate the orthonormality constraint into the variational principle via Lagrange multipliers. This follows: 
\begin{equation}\label{Eqn:min_princ}
 \mathcal{G}(\dot{\bm{U}}, \dot{\bm{Y}}) = \left\Vert \mathbb{E}\bigg[ \sum_{i=1}^r\frac{d (\bm{u}_i \bm{y}^T_i)}{dt}  - \dot{\bm{T} } \bigg]^{2}  \right\Vert_F + \sum_{i,j=1}^r \lambda_{ij}(t)  \big( \inner{\bm{u}_i}{\dot{\bm{u}}_j} - \bs{\phi}_{ij} \big),
\end{equation}
where $\lambda_{ij}(t), i,j=1, \dots, r$ are the Lagrange multipliers.  
In Appendix A, we show  that the first-order optimality condition of the above functional leads to the closed-form evolution equation of $\bm{U}$ and $\bm{Y}$ as given by:
\begin{align}
\dot{\bm{U}} &=  \underset{\perp \bm{U}}{\prod } \mathbb{E}[\dot{\bm{T}} \bm{Y}] \bm{C}^{-1}+  \bm{U}\bs{\phi} \label{eq:u_ev},\\
\dot{\bm{Y}} &= \inner{\dot{\bm{T}^T}}{ \bm{U}} + \bm{Y} \bs{\phi} \label{eq:y_ev},
\end{align}
where
$\bm{C} = \mathbb{E}[\bm{Y}^T \bm{Y}]$ is the reduced covariance matrix ($\bm{C} \in \mathbb{R}^{r\times r} $), and $\underset{\perp \bm{U}}{\prod }$ is the orthogonal projection matrix to the space spanned by the complement of $\bm{U}$, i.e. 
$\underset{\perp \bm{U}}{\prod } \bm{V} = \bm{V} - \bm{U}\bm{U}^T\Wx\bm{V} = (\bm{I}-\bm{U}\bm{U}^T\Wx)\bm{V}$.\\

Solving equations \ref{eq:u_ev} and \ref{eq:y_ev} require determination of $\bs{\phi}$. We show in Appendix B that for any skew-symmetric matrix $\bs{\phi}$ the above reduction remains equivalent for all times. Equivalent reductions can be defined as follows:  
The two reduced subspaces $\{\bm{Y} , \bm{U}  \}$ and $\{\widetilde{\bm{Y}} , \widetilde {\bm{U}}  \}$, where  $ \bm{U} , \widetilde {\bm{U}} \in \mathbb{R}^{n\times r}$ and   $\bm{Y},\widetilde{\bm{Y}}   \in \mathbb{R}^{s\times r} $, 
are equivalent if there exists a transformation matrix $\bm{R} \in \mathbb{R}^{r \times r}$ such that $\bm{U}=\widetilde{\bm{U}}\bm{R}$ and $\bm{Y}=\widetilde{\bm{Y}}\bm{R}$, where $\bm{R}$  is an orthogonal rotation matrix, i.e. $\bm{R}^T\bm{R}=\bm{I}$. Therefore, changing  $\bs{\phi}$ does not affect the subspace the both $\bm{U}$ and $\bm{Y}$ span, and it will only lead to an in-subpsace rotation.   In the rest of this paper, we take the simplest case of $\bs{\phi} \equiv  \bm{0}$. In that case the evolution of the dynamic basis reduces to:
\begin{align}
\dot{\bm{U}} &=  \underset{\perp \bm{U}}{\prod } \mathbb{E}[\dot{\bm{T}} \bm{Y}] \bm{C}^{-1} \label{eq:u_ev_DO_temp},\\
\dot{\bm{Y}} &= \inner{\dot{\bm{T}^T}}{ \bm{U}} \label{eq:y_ev_DO_temp}.
\end{align}
For the computation of the right hand side of equation \ref{eq:u_ev_DO_temp}, $\bm{U}\bm{U}^T$ may  not be computed nor stored explicitly. This can be done by observing that:
\begin{equation*}
 \underset{\perp \bm{U}}{\prod } \mathbb{E}[\dot{\bm{T}} \bm{Y}] \bm{C}^{-1} = \lp \bm{I} -   \bm{U}\bm{U}^T \Wx \rp  \dot{\bm{T}} \Wxi \bm{Y} \bm{C}^{-1} = \lp \dot{\bm{T}} \Wxi \bm{Y}  - \bm{U}\dot{\bm{T}_r} \rp \bm{C}^{-1},
\end{equation*}
where $\dot{\bm{T}_r} = \bm{U}^T \Wx  \dot{\bm{T}} \Wxi \bm{Y} $ and $\dot{\bm{T}_r} \in \mathbb{R}^{r\times r}$.
Therefore, the evolution equations for the dynamic basis and the reduced manifold are given by:
\begin{align}
\dot{\bm{U}} &=  \lp \dot{\bm{T}} \Wxi \bm{Y}  - \bm{U}\dot{\bm{T}_r} \rp \bm{C}^{-1} \label{eq:u_ev_DO},\\
\dot{\bm{Y}} &=  \dot{\bm{T}}^T \Wx    \bm{U}. \label{eq:y_ev_DO}
\end{align}

Besides the simplicity of the above evolution equations, the choice of $\bs{\phi}=\bm{0}$ also leads to an insightful interpretation of the dynamic basis reduction, namely the  update of the dynamic basis $\dot{\bm{U}}$ is always orthogonal to the space spanned by $\bm{U}$, i.e. $\dot{\bm{U}} \perp \bm{U}$, and this update is  driven by the projection of $\mathbb{E}[\dot{\bm{T}}\bm{Y}]$ onto to the orthogonal complement of $\bm{U}$ (see equation \ref{eq:u_ev_DO}). This means that if $\mathbb{E}[\dot{\bm{T}}\bm{Y}]$ is in the span of $\bm{U}$, then
\begin{equation*}
 \dot{\bm{U}} = \underset{\perp \bm{U}}{\prod } \mathbb{E}[\dot{\bm{T}} \bm{Y}] \bm{C}^{-1} = \bm{0},
 \end{equation*}
 However, when  $\mathbb{E}[\dot{\bm{T}}\bm{Y}]$ is not in the  span of $\bm{U}$, it results in $\dot{\bm{U}} \neq  \bm{0}$.
 This built-in adpativity of the dynamic basis allows $\bm{U}$ to ``follow" the data instantaneously. The  $\bm{Y}$ coefficients are updated by  projection of $\dot{\bm{T}}$ onto the space spanned by $\bm{U}$ (see equation \ref{eq:y_ev_DO}). Therefore, the evolution equation of $\bm{Y}$ captures the evolution  of $\bm{T}$ within the subspace $\bm{U}$. We also observe that the above scheme is agnostic with regard to the underlying model that generates the observation and observations may come from a system whose model is not known.

\subsection{Workflow and computational complexity}
The dynamic basis reduction is governed by time-continuous evolution equations \ref{eq:u_ev_DO} and \ref{eq:y_ev_DO}. For time-discrete data these equations must be discretized in time.  Here we summarize all the steps required for solving equations \ref{eq:u_ev_DO} and \ref{eq:y_ev_DO}:\\

\textbf{Step 1:} Given the  time-discrete realizations in the form of $\widetilde{\bm{T}}_k, k=0,1,2, \dots, K$ we compute zero-mean by subtracting the instantaneous sample mean from the data: $\bm{T}_k   = \widetilde{\bm{T}}_k - \mathbb{E}[\widetilde{\bm{T}}_k] \mathbf{1}_{1\times s}$, where $\mathbf{1}_{1\times s}$ is a row vector of size $s$ whose all elements are equal to one.\\

\textbf{Step 2:} We compute the initial condition for $\bm{Y}_0$ and $\bm{U}_0$ by performing  Karhunen-Lo\`eve decomposition of the  $\bm{T}_0$
at $t=0$. This follows:
\begin{equation*}
\bm{T}_0 = \sum_{i=1}^r \sqrt{\lambda_i} \bm{v}_i \bm{z}_i^T + \bm{E}_0,
\end{equation*}
where vectors $\bm{v}_i \in \mathbb{R}^n$ are an orthonormal basis, $\bm{z}_i \in \mathbb{R}^s$ are the stochastic coefficients, and   $\lambda_i$ are their  correspond eigenvalues and $\bm{E}_0$ is the reduction error. The initial conditions for $\bm{U}(t)$ and the stochastic coefficients $\bm{Y}(t)$ are then given by: $\bm{u}_{i_0} = \bm{v}_i$ and $\bm{y}_{i_0}= \sqrt{\lambda_i} \bm{z}_i$ for $i=1,2, \dots, r$. \\

\textbf{Step 3:} The time derivative of the observations   $\dot{\bm{T}}(t)$ can be estimated via finite difference schemes. For example the fourth-order central finite difference approximation takes the form of:
\begin{equation*}
\dot{\bm{T}}_k = \frac{-\bm{T}_{k+2} + 8\bm{T}_{k+1} - 8\bm{T}_{k-1} + \bm{T}_{k-2} }{12 \Delta t} + \mathcal{O}(\Delta t^4).
\end{equation*}
For the first two time steps one sided finite-difference schemes may be used.   \\

\textbf{Step 4:} Evolve the dynamic basis and the stochastic coefficients to the next time step $k$ using a time-integration scheme. We use fourth-order Runge-Kutta scheme. A comparison with lower-order time integration is scheme is given in the last demonstration case. \\

\textbf{Step 5:} Perform a Gram–Schmidt orthonormalization on the dynamic basis $\bm{U}$. Although the dynamic basis remains orthonormal by construction, the numerical orthonormality of $\bm{U}$ may not be preserved due to the nonlinearity of the evolution equations \ref{eq:u_ev_DO} and \ref{eq:y_ev_DO}. An  enforcement of orthonormality of $\bm{U}$ resolves this issue. Now steps 3-5 can be repeated  in a time-loop till  the completion of all time steps.   

A quick inspection of the evolution equations \ref{eq:u_ev_DO}-\ref{eq:y_ev_DO} reveals two favorable properties with regard to the computational cost and memory storage requirement of  the dynamic basis reduction: (1)  
The computational complexity of solving the evolution equations \ref{eq:u_ev_DO}-\ref{eq:y_ev_DO} is $\mathcal{O}(rns)$, i.e. linear computational complexity with respect to the reduction size $r$, the size of the observations in state space $n$ and the sample size $s$.  We note that  given that $r$ is often a small number, i.e.  $r << n$ and $r << s$,  the cost of inverting the covariance matrix $\mathbf{C}$
is often negligible. (2) The memory storage requirement is also linear with respect to $n$, $r$ and $s$. This is because at each time instant only computation of $\dot{\bm{T}}$ is required, which can be computed from few snapshots of $\bm{T}$. This allows for loading  matrix $\bm{T}$ in small time batches. The linear computational and memory complexity of the dynamic basis reduction may be exploited for applications with very large data sets.   

\subsection{Mode ranking within the subspace}
In this section we present a ranking criterion within the dynamic subspace $\bm{U}$. This can be accomplished by using the reduced covaraince matrix $\bm{C}$, which, in general, is a full matrix and this implies that the $\bm{Y}$ coefficients are in general correlated. The uncorrelated $\bm{Y}$ coefficients may be computed via the eigen-decomposition of the covariance matrix $\bm{C}$. To this end, let
\begin{equation}
\bm{C} \bs{\psi} = \bs{\psi} \bm{\Lambda},
\end{equation}
be the eigenvalue decomposition of matrix $\bm{C}, $ where $\bs{\psi} \in \mathbb{R}^{r\times r}$ is matrix of eigenvectors of $\bf{C}$ and $\bs{\Lambda}=\mbox{diag}(\lambda_1, \lambda_2, \dots, \lambda_r)$ is a diagonal matrix of the eigenvalues of matrix $\bm{C}$. The uncorrelated $\bm{Y}$ coefficients and the dynamic basis associated with those uncorrelated coefficients, denoted by $\hat{\bm{Y}}$ and $\hat{\bm{U}}$, are obtained by an in-subspace rotation:
\begin{equation*}
    \hat{\bm{Y}} = \bm{Y}\bs{\psi}, \quad \quad \hat{\bm{U}} = \bm{U}\bs{\psi}.
\end{equation*}
It is easy to verify that  $\hat{\bm{Y}}$ coefficients  are uncorrelated, i.e. $\mathbb{E}[\hat{\bm{Y}}\hat{\bm{Y}}^T]=\bs{\Lambda}$. We adopt a variance-based ranking for  $\{\hat{\bm{Y}}, \hat{\bm{U}}\}$.  Since matrix $\bm{C}$ is  positive semidefinite, without loss of generality, let $\lambda_1 \geq \lambda_2 \geq ... \geq  \lambda_r \geq 0$, where $\lambda_i$ is the variance of $\hat{\bm{y}}_i$. 

 Since  $\{\hat{\bm{Y}}, \hat{\bm{U}}\}$  are obtained by an in-subspace rotation of $\{\bm{Y}, \bm{U}\}$, they are equivalent representations of the low rank dynamics, i.e. $\hat{\bm{U}}\hat{\bm{Y}}^T=\bm{U}\bs{\psi} \bs{\psi}^T\bm{Y}^T = \bm{U}\bm{Y}^T$, where we have used $\bs{\psi} \bs{\psi}^T = \bm{I}$, since eigenvectors of a symmetric matrix ($\bm{C}$) are orthonormal. In all the demonstration cases in this paper, we use the  ranked uncorrelated representation of the dynamic basis for the analysis of the results.   For the sake of simplicity in the notation, we use $\{\bm{Y}, \bm{U}\}$ to refer to the ranked uncorrelated representation instead of  $\{\hat{\bm{Y}}, \hat{\bm{U}}\}$, noting that both of these representations are equivalent.  

\section{Demonstration cases}
\subsection{Stochastic advection equation}
As the first example, we apply dynamic basis reduction to  a one-dimensional stochastic advection equation.  As we will show in this section, the stochastic advection equation chosen here is exactly two-dimensional, when expressed versus the dynamic basis. It also has an exact analytical  solution for the dynamic basis $\bm{U}$ and the stochastic coefficients $\bm{Y}$, which enables an appropriate verification of the dynamic basis algorithm.  The stochastic advection equation is given by:
\begin{equation}\label{eq:stoch_adv}
\pfrac{u}{t} + V(\xi) \pfrac{u}{x} = 0,  \quad x \in[-1,1], \quad t \in[0,T_f],
\end{equation} 
 subject to periodic boundary conditions at the ends of the spatial domain $x\pm 1$ and the initial condition  of $u_0(x;\xi)$. In the above equation $T_f$ is the final time.  The solution to the above equation is a stochastic random field in the form of $u(x,t;\xi)$. The stochastic advection speed $V(\xi)$ and the deterministic initial condition  are chosen to be
\begin{equation*}
V(\xi)=\ol{v} + \sigma \xi, \quad \quad u_0(x;\xi) = \sin \pi x, 
\end{equation*} 
where $\xi$ is a one-dimensional uniform random variable defined on $[-1,1]$ and $\sigma$ is a constant and
$\ol{v}$ is the mean of the transport velocity. The exact solution of equation \ref{eq:stoch_adv} is
\begin{equation*}
u(x,t;\xi) =  \sin \pi \big(x- (\ol{v}+\sigma \xi)t \big), 
\end{equation*} 
and its mean and variance are:
\begin{align*}
& \mathbb{E}[u(x,t;\xi)] =  \frac{1}{\sigma \pi t} \sin \pi (x- \ol{v}t )
\sin \sigma \pi t,\\
& \mathbb{E}[u(x,t;\xi)^2] - \mathbb{E}[u(x,t;\xi)]^2 =  \frac{1}{2} 
                                                     -  \frac{1}{4 \sigma
\pi t} \cos 2\pi (x- \ol{v}t ) \sin 2 \sigma \pi t
                                                     -  \frac{1}{\sigma^2
\pi^2 t^2} \sin^2 \pi (x- \ol{v}t ) \sin^2 \sigma \pi t, 
\end{align*} 
The expectation operator is given by $\mathbb{E}[u(x,t;\xi)] = \displaystyle \int_{-1}^1 u(x,t;\xi) \rho(\xi) d\xi  $, where $\rho(\xi) = 1/2$ is the probability density function. The total variance is equal to:
\begin{equation*}
\Sigma(t) = \int_{-1}^1 \big( \mathbb{E}[u(x,t;\xi)^2] - \mathbb{E}[u(x,t;\xi)]^2 \big) dx = 1 - \frac{\sin^2 \sigma \pi t}{(\sigma \pi t)^2}.
\end{equation*}
The mean-subtracted stochastic field can be expressed with the following two dimensional expansion:
\begin{equation}\label{eq:adv_exp}
u(x,t;\xi) - \mathbb{E}[u(x,t;\xi)] = u_1 y_1 + u_2 y_2, 
\end{equation} 
where 
\begin{align*}
u_1(x,t) &= \sin \pi(x- \ol{v}t ),          &y_1(\xi,t) = \cos\pi \sigma\xi t - \frac{1}{\sigma \pi t} \sin \pi \sigma t, \\
u_2(x,t) &= \cos \pi(x- \ol{v}t ),          &y_2(\xi,t) =  \sin \pi \sigma\xi t.
\end{align*}
It is easy to verify that the  two-dimensional expansion of the above random field satisfies the dynamic basis orthonormality condition in the continuous form, i.e.
\begin{equation*}
 \int_{-1}^1 u_i(x,t) u_j(x,t) dx = \delta_{ij}, \qquad i,j=1,2.
\end{equation*} 
The expansion of the stochastic field as expressed above in Equation \ref{eq:adv_exp} is in bi-orthogonal form since 
$\mathbb{E}[y_i y_j]=0, i\neq j$. 

The details of solving the evolution equations of the dynamic basis $\bm{U}(t)$ and the stochastic coefficients $\bm{Y}(t)$ are as follows. The matrix-valued random observations $\bm{T}(t)$ are generated by solving the stochastic advection equation \ref{eq:stoch_adv} at $s=64$ collocation points. Theses sampling points are chosen to be the roots of the Legendre polynomial of order 64 defined on the interval of $\xi \in [-1,1]$. The 64 realizations of equation \ref{eq:stoch_adv}  are generated using a deterministic solver that uses Fourier modes of order $n=128$ for space discretization and fourth-order Runge-Kutta with $\Delta t= 0.001$ for time-advancement. The final time is set to  $T_f=10$ for all realizations.  

\begin{figure}
\centering
\subfigure[]{
\includegraphics[width=.47\textwidth]{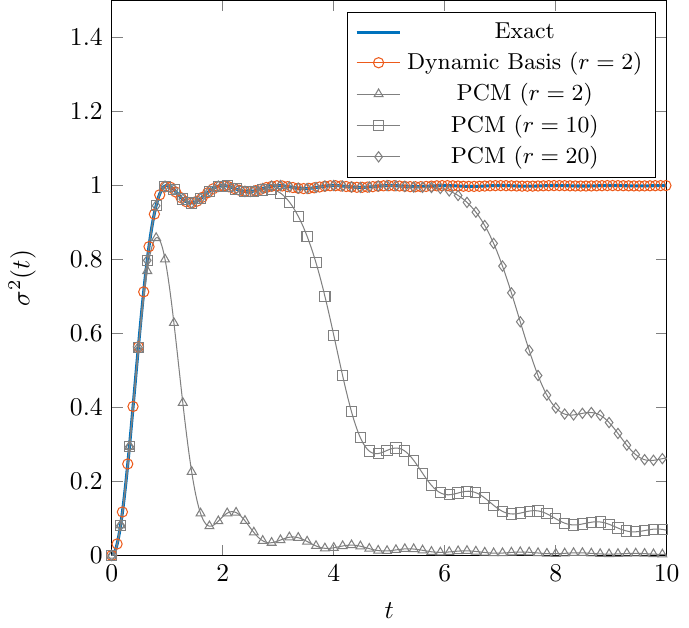}
\label{Fig:Adv-Total-Var}
}
\subfigure[]{
\includegraphics[width=.47\textwidth]{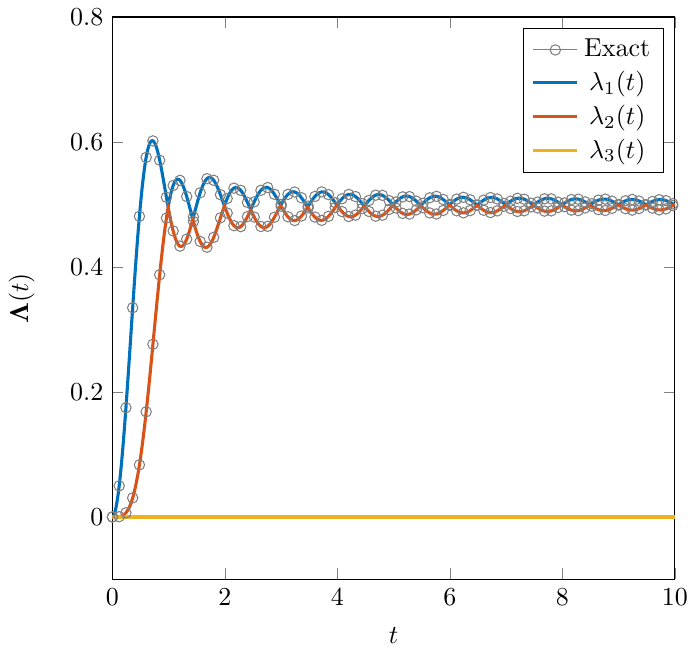}
\label{Fig:Adv-Eig-r_3}
}
\caption{Stochastic advection equation: (a) Comparison of the total variance of the dynamic basis reduction and polynomial chaos method (PCM) with the exact solution. (b) Comparison of the eigenvalues of the covariance matrix obtained from dynamic basis reduction with those of the exact solution. The exact solution is exactly two-dimensional when expressed in terms of the dynamic basis, and the variance picked up by the third mode is zero (order of round off error).   }
\end{figure}

We first performed dimension reduction of size $r=2$. In Figure \ref{Fig:Adv-Total-Var}, the total variance $\Sigma(t)$ obtained from the reduction is compared against  the exact variance. The total variance  of the low-rank approximation is obtained by $\Sigma(t) = \sum_{i=1}^r \lambda_i(t)$ where $\lambda_i(t)$ are the eigenvalues of the covariance matrix $\bm{C}$. We also solved the stochastic advection equation \ref{eq:stoch_adv} using probabilistic collocation method (PCM), in which the stochastic solution is approximated as:
\begin{equation*}
u(x,t;\xi) - \mathbb{E}[u(x,t;\xi)] = \sum_{i=1}^{r} \phi_i(\xi) v_i(x,t) + e(x,t;\xi),
\end{equation*}
where $\phi_i(\xi)$ are Legendre-chaos polynomials of order $i$, $v_i(x,t)$ are space-time coefficients and $e(x,t;\xi)$ is the expansion error.  In Figure \ref{Fig:Adv-Total-Var}, the total variance obtained from different reduction sizes $r$ of the PCM are also shown where it can be seen that as the size of PCM expansion increases the  PCM variance follows the exact solution for a longer time, but eventually the PCM error increases and even a 20-dimensional  PCM reduction under-performs the dynamic basis reduction with $r=2$. For the long-term error analysis for this problem using the polynomial chaos expansion, the reader may refer to reference \cite{WK06}.  We note that in the PCM reduction the stochastic basis is fixed \emph{ a priori}, i.e. the polynomial chaos $ \phi_i(\xi)$ while in the dynamic basis reduction the stochastic coefficients $\bm{y}_i(\xi,t)$ are time-dependent and intuitively this enables the stochastic coefficient adapt instantaneously to the data   according to the variational principle  \ref{Eqn:min_princ}.

\begin{figure}
\centering
\includegraphics[width=.68\textwidth]{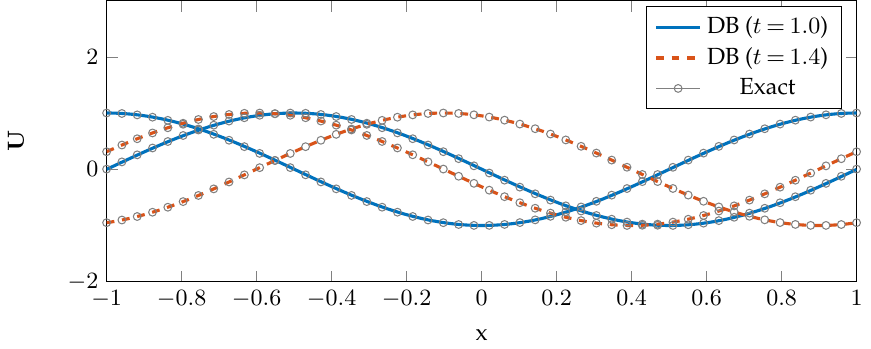}
\includegraphics[width=.28\textwidth]{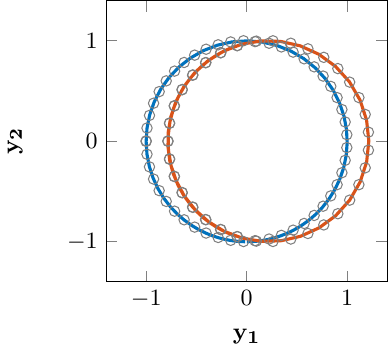}
\caption{Reduction of the solution of the stochastic advection equation using dynamics basis: (a) Dynamic basis; (b) Stochastic coefficient (phase portrait).}
\label{Fig:Adv_U_Y}
\end{figure}

In Figure \ref{Fig:Adv-Eig-r_3} the covariance eigenvalues of the dynamic basis reduction with $r=3$ are compared against the exact solution. 
Since the solution is two-dimensional when expressed versus the dynamic basis, the third mode must remain inactive. As shown in Figure   \ref{Fig:Adv-Eig-r_3}, this is confirmed numerically as the third eigenvalue  does not pick up energy, while the first two eigenvalues are in excellent agreement with the exact solution. In general small eigenvalues of the covariance matrix may cause numerical instability  due to the  computation of $\bm{C}^{-1}$ in equation \ref{eq:u_ev_DO}. The issue of singularity of matrix $\bm{C}$ can be resolved by computing the pseudo-inverse of the covaraince matrix. See reference \cite{Babaee:2017aa} for similar treatment of singularity of the covaraince matrix.  In  Figure \ref{Fig:Adv_U_Y}, the dynamic basis and the stochastic coefficients are compared against the analytical values for two different times. In both cases excellent agreement is observed.

\subsection{Transient Instability:  Kuramoto-Sivashinsky }
\begin{figure}
\centering
\subfigure[$u(x,t)$]{
\includegraphics[width=.29\textwidth]{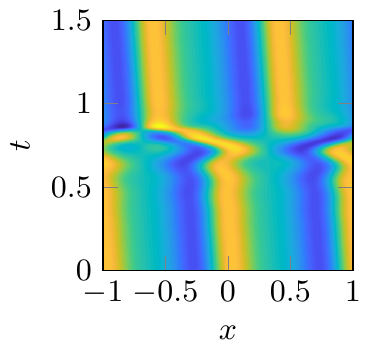}
\label{Fig:KS_u}
}
\subfigure[$K(x,x')$, $l_c=2.5$]{
\includegraphics[width=.31\textwidth]{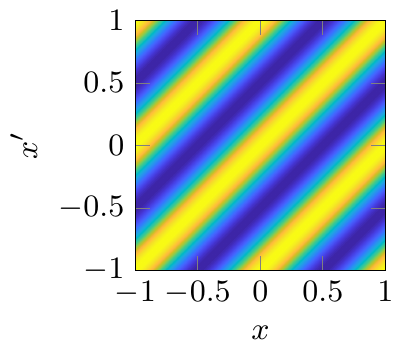}
\label{Fig:KS_lc=25e-1}
}
\subfigure[$K(x,x')$, $l_c=0.05$]{
\includegraphics[width=.31\textwidth]{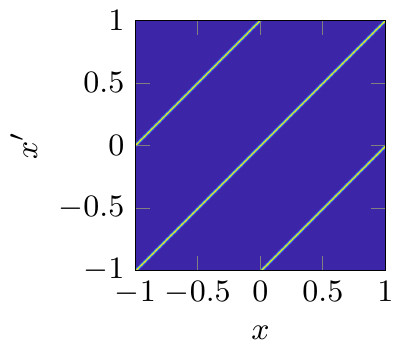}
\label{Fig:KS_lc=5e-2}
}
\caption{Kuramato-Sivashinsky:(a) Solution $u(x,t)$ shows a sudden transition;  covariance kernel of the random perturbation with (b) $l_c=2.5$, and (c) $l_c=0.05$.}
\end{figure}

One of the main utilities of a low-rank approximation of high-dimensional dynamical systems  is to provide  interpretable explanation of complex spatio-temporal dynamics.  In this section, we demonstrate the application of the dynamic basis reduction to the one-dimensional stochastic Kuramoto-Sivashinsky equation. This equation has been widely studied for chaotic spatio-temporal dynamics in wide-ranging applications including reaction-diffusion systems \cite{KT76}, 
fluctuations in fluid films on inclines \cite{SM80} and  instabilities in laminar flame fronts \cite{S77}.  The Kuramoto-Sivashinsky equation exhibits \emph{intrinsic instabilities} and as a result the stochastic perturbations may undergo a significant growth. The spontaneity  of the instability makes it a challenging problem for any reduction technique as the reduced subspace during the transition could be very different from that of the equilibrium states. We apply the dynamic basis reduction to realizations of the Kuramoto-Sivashinsky system  to demonstrate how the $\bm{U}(t)$   ``follows'' the dynamics as the observations undergo a transient instability. This transient instability is  visualized in Figure \ref{Fig:KS_u} where a sudden transition occurs near $t=0.75$.  The solution shown in Figure \ref{Fig:KS_u} is one realization of the Kuramoto-Sivashinsky equation with the following setup. We consider Kuramoto-Sivashinsky equation
\begin{align}\label{eqn:KS}
& \pfrac{u}{t} = u\pfrac{u}{x} - \pfrac{^2 u}{x^2}-\epsilon \pfrac{^4 u}{x^4},
&  x \in[-1,1], \quad t \in[0,T_f],
\end{align}
with $\epsilon=0.01$ and subject to random initial condition
\begin{align}
& u(x,0;\xi)   = u_b(x) + u'(x;\xi), & x \in[-1,1].
\end{align}
In the above equation $u_b(x)$ is chosen to be the terminal solution of equation \ref{eqn:KS} after long-time integration.
The random perturbation of the initial condition, denoted by $u'(x;\xi)$, is set to have  the periodic covariance given by:
\begin{equation*}
K(x,x') = \sigma^2 \exp{\lp -\frac{2\sin^2(\pi(x-x'))}{l_c^2} \rp}.
\end{equation*}
where $\sigma$ is the standard deviation and $l_c$ is the correlation length.
 We perform Karhunen-Lo\`eve decomposition 
\begin{equation*}
 \int_{-1}^1 K(x,x')\phi_i(x')dx' = \lambda_i  \phi_i(x), 
\end{equation*}
where $\lambda_i$ and $\phi_i(x)$ are the $i$th eigenvalue and eigenfunction of the covariance kernel.  
The random perturbation is then approximated as:
\begin{equation*}
u'(x;\xi) = \sum_{i=1}^d \sqrt{\lambda_i}\xi_i \phi_i(x) + \epsilon,
\end{equation*}
where  $\xi_i \in [0,1]$ are independent uniform random numbers with unit variance. 
The dimension $d$ is chosen such that 99\% of the energy is captured by the Karhunen-Lo\`eve decomposition.  Decreasing $l_c$ increases
the number of dimensions required to capture the energy threshold.  
We consider two cases of low-dimensional and high-dimensional random perturbations with $l_c=2.5$ and $l_c=0.05$, respectively.
The covariance kernels for these two correlation lengths are shown in Figures \ref{Fig:KS_lc=25e-1}-\ref{Fig:KS_lc=5e-2}.

 \textbf{Three-dimensional random perturbation:} We consider the correlation length of $l_c=2.5$ and  $\sigma=0.1$. 
The three-dimensional ($d=3$) expansion of the KL decomposition captures more than $99\%$ of the energy of the covariance kernel. We sample the Kuramoto-Sivashinsky equation at probabilistic $Q=5$ collocation points  in each $\xi_i$ direction, and therefore,  the total number of samples is $s=5^3$. Correspondingly, the sample weights  $\bm{w_{\xi}}$ are the probabilistic collocation weights.  We use \emph{chebfun} \cite{chebfun}  as a blackbox solver to solve the Kuramoto-Sivashinsky equation \ref{eqn:KS} for $s=5^3$ realizations of the random initial condition. The solution is post-processed  with fixed Chebyshev polynomial order $n=256$. Matrix $\bm{T}$ is sampled at discrete times with constant time interval of $dt=0.001$ from $t=0$ to $T_f=1.2$. This time interval encompasses the transition period that occurs $0.6<t<1.0$.

\begin{figure}[t]
\centering
\subfigure[Eigenvalues of $\bm{C}$]{
\includegraphics[width=.45\textwidth]{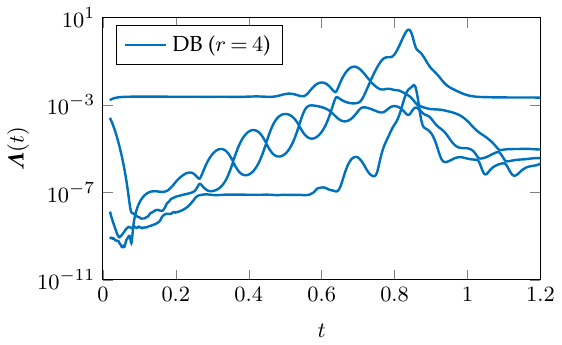}
\label{Fig:KS_eig_d=3}
}
\subfigure[Dynamics Basis 1]{
\includegraphics[width=.45\textwidth]{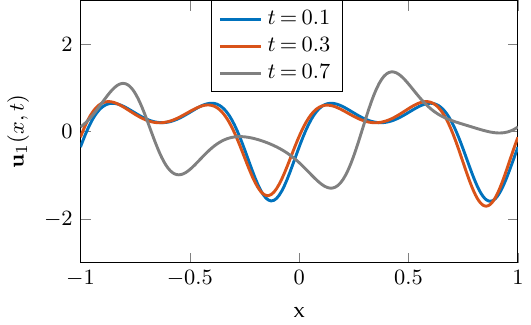}
\label{Fig:KS-mode1}
}
\subfigure[Dynamics Basis 2]{
\includegraphics[width=.45\textwidth]{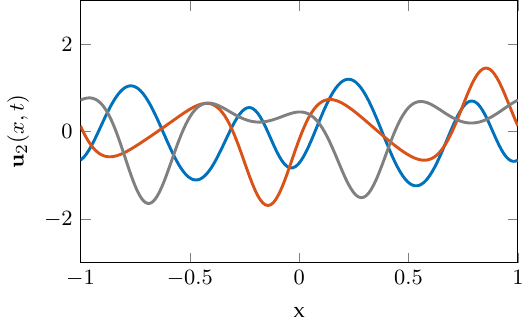}
\label{Fig:KS-mode2}
}
\subfigure[Dynamics Basis 3]{
\includegraphics[width=.45\textwidth]{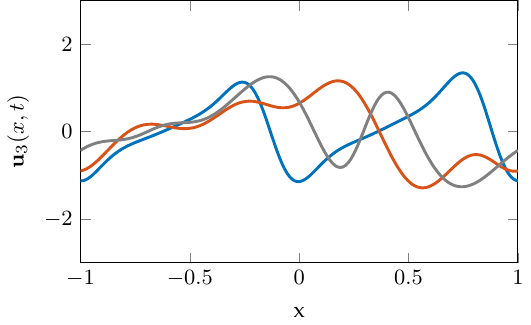}
\label{Fig:KS-mode3}
}
\caption{Dynamic basis reduction for Kuramoto-Sivashinsky equation with the reduction  size $r=4$: (a) eigenvalues of reduced covaraince matrix.   The first three modes at different time instants: (b) $\bm{u}_1$, (c) $\bm{u}_2$  and (d) $\bm{u}_3$. }
\end{figure}

We perform dynamic basis reduction with size $r=4$. In Figure  \ref{Fig:KS_eig_d=3}, the eigenvalues of the reduced covariance matrix versus time are shown.  The eigenvalues of the covariance matrix represent variance or the energy associated with  the dynamic basis. The first three dynamic basis associated with three largest eigenvalues of the covaraince matrix are shown in Figures \ref{Fig:KS-mode1}-\ref{Fig:KS-mode3}.
We first observe that the energy of the first  mode remains  roughly constant for $0<t<0.5$.   However, the second and third modes undergo four order of magnitude growth during the same period.  Correspondingly, the first dynamic basis remain nearly invariant for $0<t<0.5$ (see Figure \ref{Fig:KS-mode1}). However, $\bm{u}_2$ and $\bm{u}_3$ are evolving during the same time interval (see Figures \ref{Fig:KS-mode2} and \ref{Fig:KS-mode3}). The second and third modes have near eigenvalue crossings periodically in this interval which causes very rapid change of the modes when their eigenvalues become close.  The significant growth of the energy of the second and third mode reveals the transient unstable directions of the dynamical system. However, these two modes are ``buried'' underneath the most energetic mode, e.g. the first mode. We note that detecting the unstable directions is crucial for control of these systems \cite{BMS18}. Detecting these directions can be utilized in building precursors of sudden transitions or rare events in transient systems. 

\begin{figure}
\centering
\subfigure[]{
\includegraphics[width=0.31\textwidth]{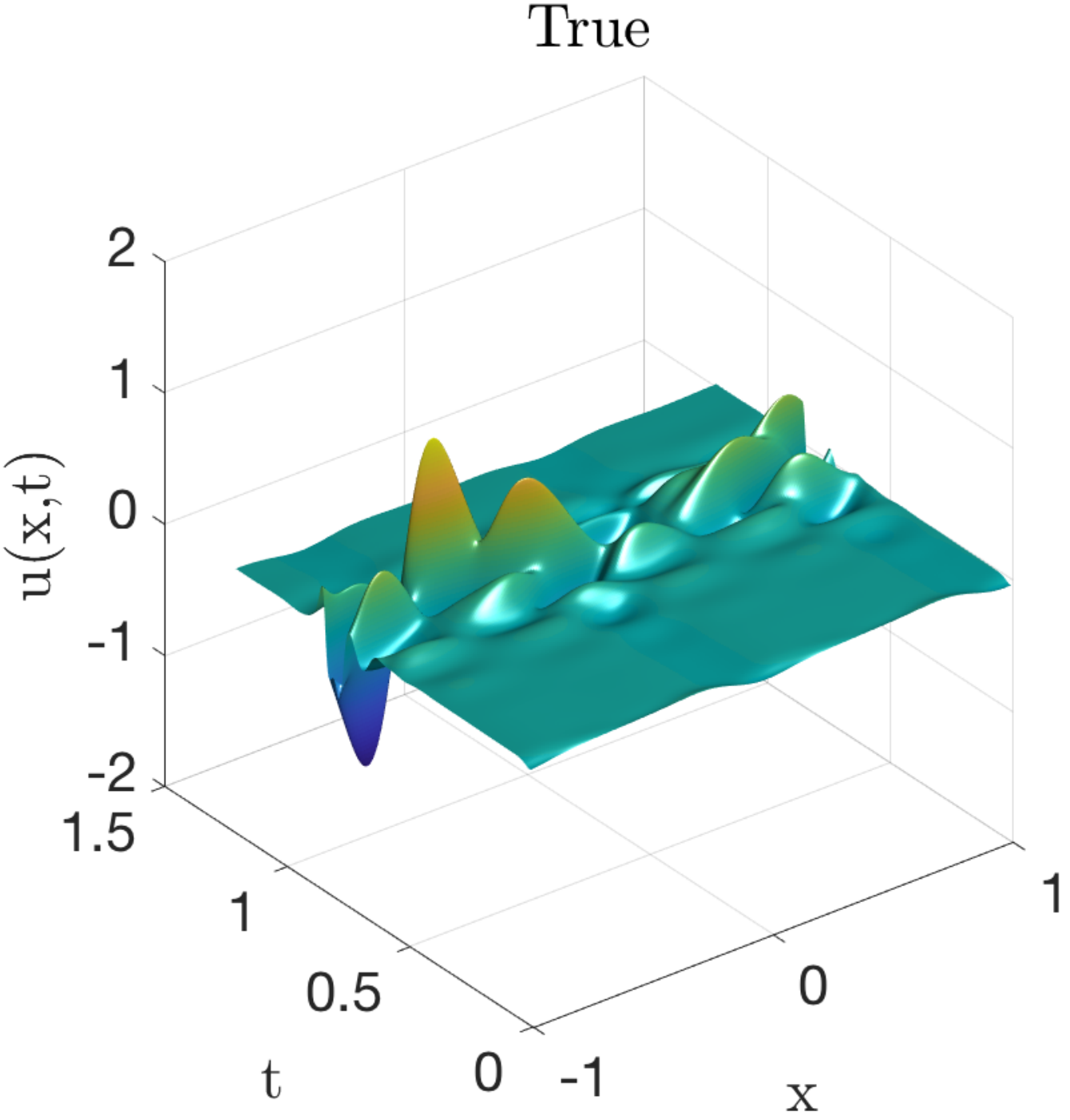}
}
\subfigure[]{
\includegraphics[width=0.31\textwidth]{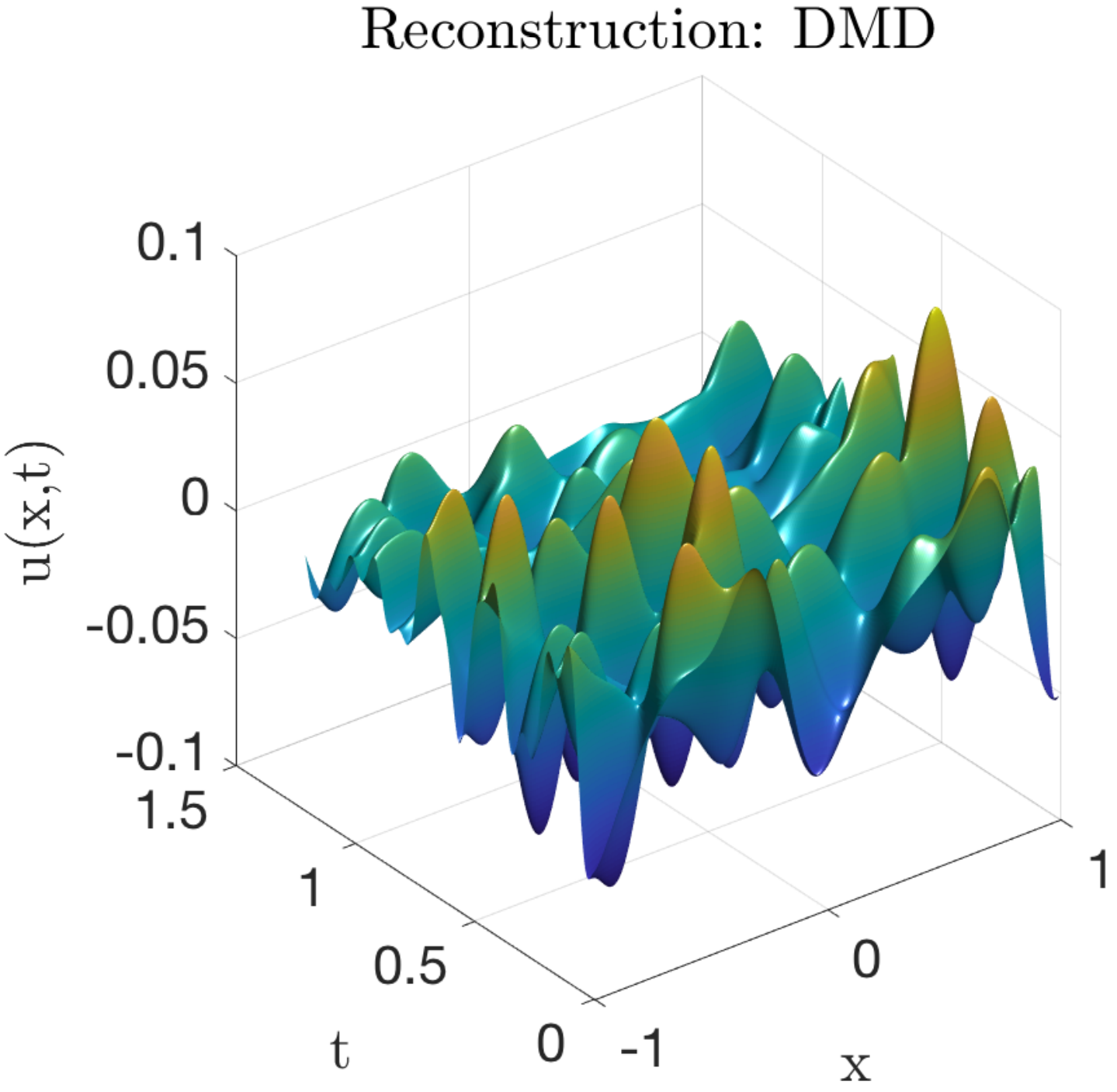}
}
\subfigure[]{
\includegraphics[width=0.31\textwidth]{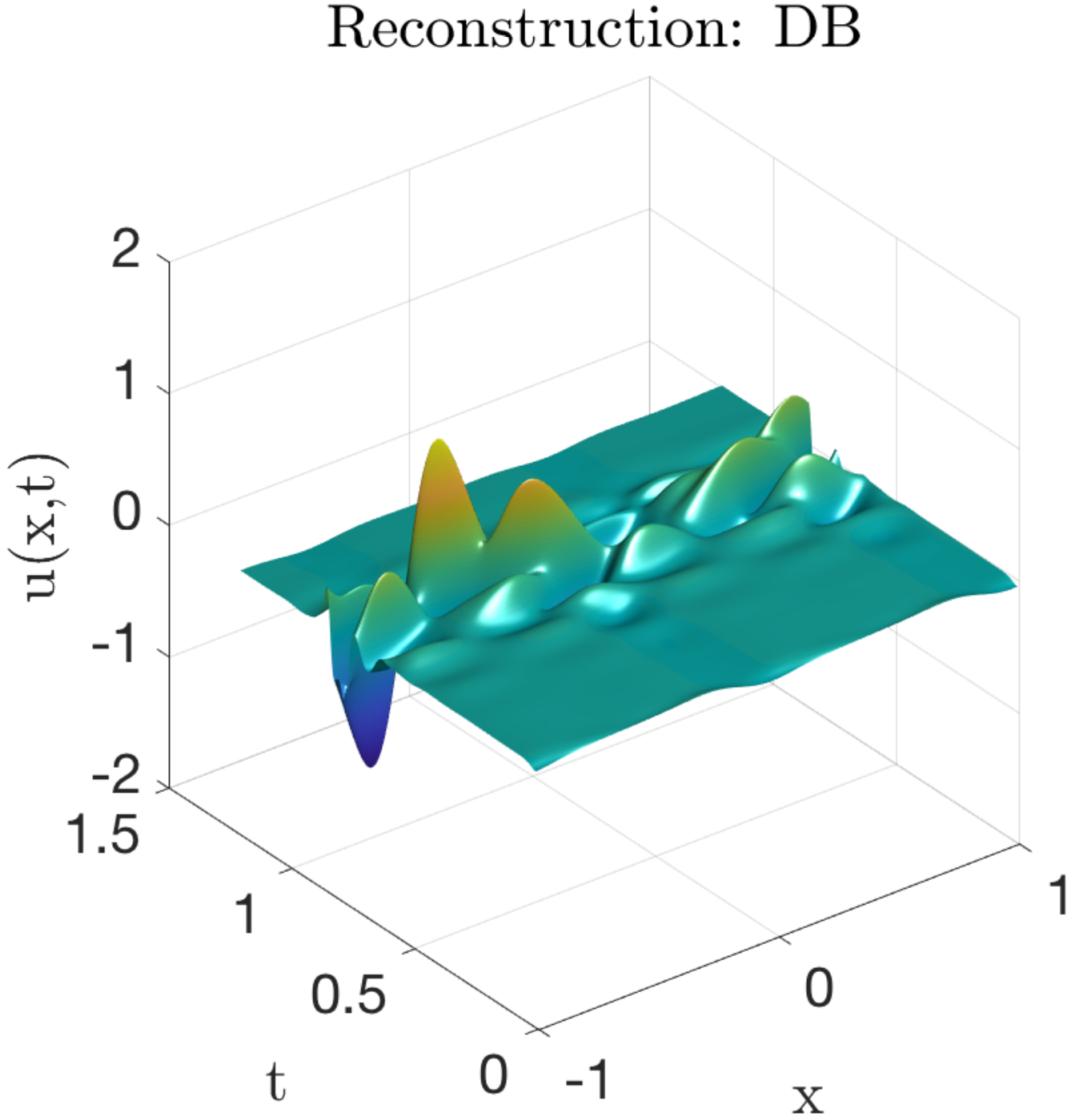}
}
\caption{Reconstruction of Kuramoto-Sivashinsky sample with: (a) True solution; (b) Reconstruction with dynamic mode decomposition (DMD); (c) reconstruction with dynamic basis (DB). }
\label{Fig:KS-DMD}
\end{figure}

We also computed DMD  modes to  compare against the dynamic basis reduction. DMD has been very successful as a diagnostic tool to characterize complex dynamical systems. However, it is difficult for DMD to handle transient phenomena.    To this end, we computed the DMD modes for a deterministic realization of Kuramoto-Sivashinsky equation. The realization is generated  by solving the Kuramoto-Sivashinsky equation for the collocation point $\bs{\xi}_0=(-0.91,-0.54,0.91)$. This collocation point is chosen arbitrarily and any other collocation point would result in the same qualitative observations that we make in this section. To compute the DMD modes we first performed  singular value decomposition (SVD) of the mean-subtracted snapshots of the Kuramoto-Sivashinsky equation in the span  of $t \in[0, T_f]$ for $\bs{\xi}_0=(-0.91,-0.54,0.91)$. The SVD truncation threshold $e=\sum_{i=1}^n \sigma_i/\sum_{i=1}^N \sigma_i$ was  set to $99\%$, where $\sigma_i$ are the singular values, and $N$ is the total number of snapshots.   The  $99\%$ threshold resulted in $n=11$ low-rank approximation. We then computed the DMD modes $\boldsymbol{\Phi}$ according to procedure explained in reference  \cite{S10}. The DMD-reconstructed data was then computed  using: $\bm{u}(t) \simeq \sum_{i=1}^n\boldsymbol{\Phi}_k \exp(\omega_k t) b_k $, where $b_k$ is the initial amplitude of each mode, $\boldsymbol{\Phi}_k$'s are the DMD eigenvectors  and $\omega_k$'s are the eigenvalues associated to DMD modes.  In Figure \ref{Fig:KS-DMD}, we compare the reconstruction of the DMD  with $n=11$ and DB with $r=4$. The DB reconstruction obtained by: $\bm{u}(t) \simeq  \sum_{i=1}^r \bm{Y}_{s_0,i}(t)\bm{u}_i(t)$, where $\bm{y}_i(t,\bs{\xi}_0)$ is the stochastic coefficients associated with $\bs{\xi}_0=(-0.91,-0.54,0.91)$, i.e. a fixed row ($s$) of matrix $\bm{Y}$. It is clear that DMD fails to capture the transient instability, while DB reconstruction looks qualitatively accurate. 

We perform a quantitative comparison between DMD and DB. We rank the $n=11$ DMD modes based on their growth rates from increasing to decreasing, and we perform DMD reconstruction $\bm{u}(t) \simeq \sum_{i=1}^r\boldsymbol{\Phi}_k \exp(\omega_k t) b_k $ for  $r=2,4,6$ and 8.  We also perform DB for the same reduction sizes. We then compute the reconstruction error: $\epsilon^2= \sum_{k=1}^K \|\bm{u}_{true}(t_k) - \bm{u}(t_k) \|^2/K $ where $K$ is the total number of snapshots and $\bm{u}(t_k)$ is the reconstructed solution based on DMD or DB. In Table \ref{tab:KS-DMD}, the reconstruction errors $\epsilon$ for DB and DMD are shown. DMD  is incapable of converging to the true solution as the size of DMD modes increases, however,  DB converges to the true solution as the reductions size increases. The limitation of DMD for capturing transient dynamics  has been recognized before; see for example \cite{KBBJ16}, Chapter 1, Section 1.5.2. \\

\begin{table}[]
    \centering
    \begin{tabular}{|c|c| c| c | c|}
    \hline
         r   & 2 & 4 & 6 & 8\\
        \hline
        Error (DMD) & 9.01$\times 10^{-2}$  & 9.027$\times 10^{-2}$ & 9.18$\times 10^{-2}$ & 9.24$\times 10^{-2}$\\
        \hline
        Error (DB) & 2.10$\times 10^{-1}$  & 6.92$\times 10^{-3}$ & $1.88\times 10^{-3}$ & 8.92$\times 10^{-4}$  \\
        \hline
    \end{tabular}
    \caption{Reconstruction error of Kuramoto-Sivashinsky equation: comparison between dynamic mode decomposition (DMD) and dynamic basis (DB). }
    \label{tab:KS-DMD}
\end{table}

\begin{figure}
\centering
\includegraphics[width=.48\textwidth]{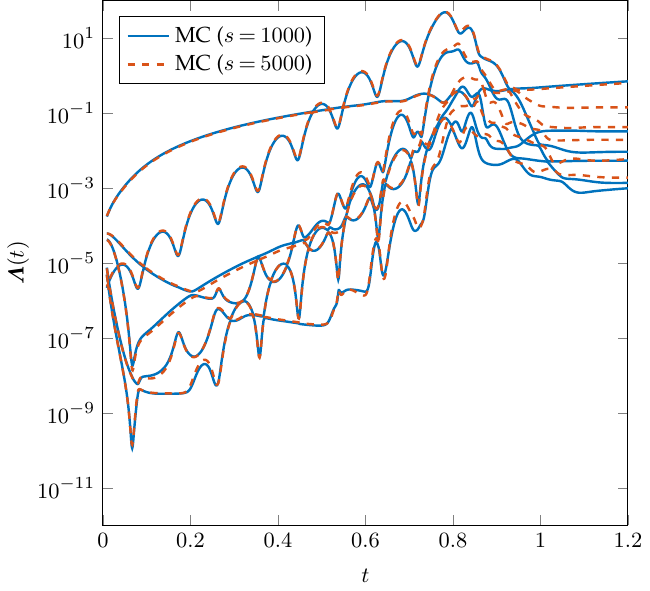}
\caption{Dynamic basis reduction of Kuramoto-Sivashinsky equation subject to one-hundred random initial conditions. Two sample sizes of $s=1000$ and $s=5000$ are shown. The reduction size is $r=6$. }
\label{Fig:KS_eig_d=100}
\end{figure}

 \textbf{One-hundred-dimensional random perturbation:}
In this section we present the results of dynamic basis reduction of   one-hundred dimensional random initial condition with correlation length of $l_c=0.05$ (see Figure \ref{Fig:KS_lc=5e-2}).  To capture the small scales of  the random initial condition the state space is represented with Chebyshev polynomial of order $n=1024$. We draw Monte Carlo samples of the initial condition with sample sizes of $s=1000$ and $s=5000$. We note these sample size are quite small given the high dimensionality of the  random space. However, as we demonstrate here the dynamic basis reduction extracts the low-dimensional structure of the observables despite the severely low number of samples. We choose reduction size of $r=6$.  In Figure \ref{Fig:KS_eig_d=100} the eigenvalues of the covariance matrix for the samples sizes are shown.  As it is clear, the dynamic basis reduction with sample sizes of 1000 and 5000 are very close before and during the transition. This implies that despite the high dimensionality of the random space, the transient response has a low-dimensional structure --- albeit a highly time-dependent one --- that is effectively captured with the proposed reduction with a severely sparse samples.  

\subsection{Transient Vertical Jet}\label{Sec:TVJ}
\begin{figure}
\centering
\subfigure[]{
\includegraphics[width=.235\textwidth]{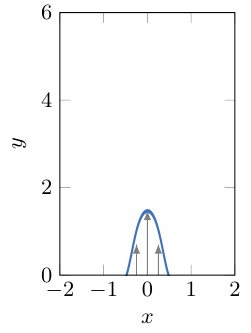}
\label{Fig:Jet_POD_Schem}
}
\subfigure[]{
\includegraphics[width=.6\textwidth]{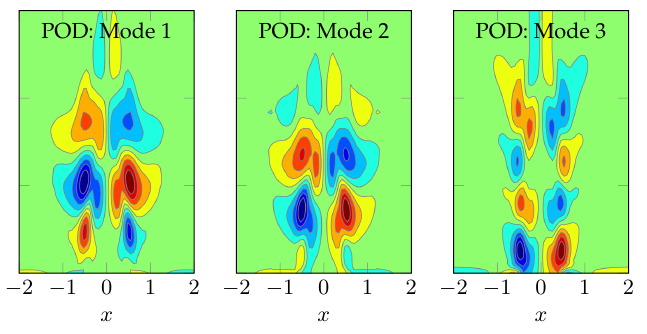}
\label{Fig:Jet_POD_Modes}
}
\caption{Transient vertical jet in quiescent flow. (a) Schematic of the problem. (b) The three most energetic POD modes.  }
\label{Fig:Jet_POD}
\end{figure}

In this section we perform dynamic reduction for a two-dimensional transient vertical jet  flushed into quiescent background flow. The jet is driven by random boundary conditions as described in this section. The schematic of the problem is shown in Figure \ref{Fig:Jet_POD_Schem}.
To generate the samples, we performed direct numerical simulation (DNS) of the two-dimensional incompressible Navier-Stokes equations:
\begin{align*}
 \pfrac{\bm{u}}{t} + (\bm{u}\cdot \nabla) \bm{u}  &= -\nabla p + \frac{1}{Re}\nabla^2\bm{u},\\
 \nabla \cdot \bm{u} &= 0,
\end{align*}
where $\bm{u}(\bm{x},t;\bxi)$ is the random velocity field  with two components of $\bm{u}=(u,v)$ and the spatial coordinates of $\bm{x}=(x,y)$, and  $p(\bm{x},t;\bxi)$ is the random pressure field  and $Re=500$ is the Reynolds number. We consider a stochastic jet velocity given by:
 \begin{equation}
 u(x,t;\bxi) = \lp 1+u'(x,\bxi) \rp  \lp 1-(x/r)^2 \rp \exp{(-.823 (x/r)^4)}.
 \end{equation}
 The mean of the prescribed jet boundary condition is similar to previous studies \cite{HessamMasters,bagheri2009global}. The random component of the boundary condition is represented by a uniform random field with the covariance kernel specified by:
\begin{equation*}
K(x,x') = \sigma^2 \exp{\lp -\frac{(x-x')^2}{l_c^2} \rp},
\end{equation*}
where $l_c=0.5$ and $\sigma=0.01$.  The dimension $d=2$ is chosen such that 99\% of the energy is captured by the Karhunen-Lo\`eve decomposition.
We performed DNS using  spectral/hp element method with $Ne=1604$ and spectral polynomial order of $p=5$. With these settings, the size of the discrete phase space is roughly $n=115,000$ when accounting for the two components of the velocity vector field. The time integration was performed using  third-order stiffly stable scheme \cite{KS05} with the time intervals of $\Delta t=5\times 10^{-3}$.  We sample the DNS at PCM points with   $Q=6$  quadrature points in each random dimension. Therefore, the total number of samples are $s=Q^2=36$.  

\begin{figure}
\centering
\includegraphics[width=1.\textwidth]{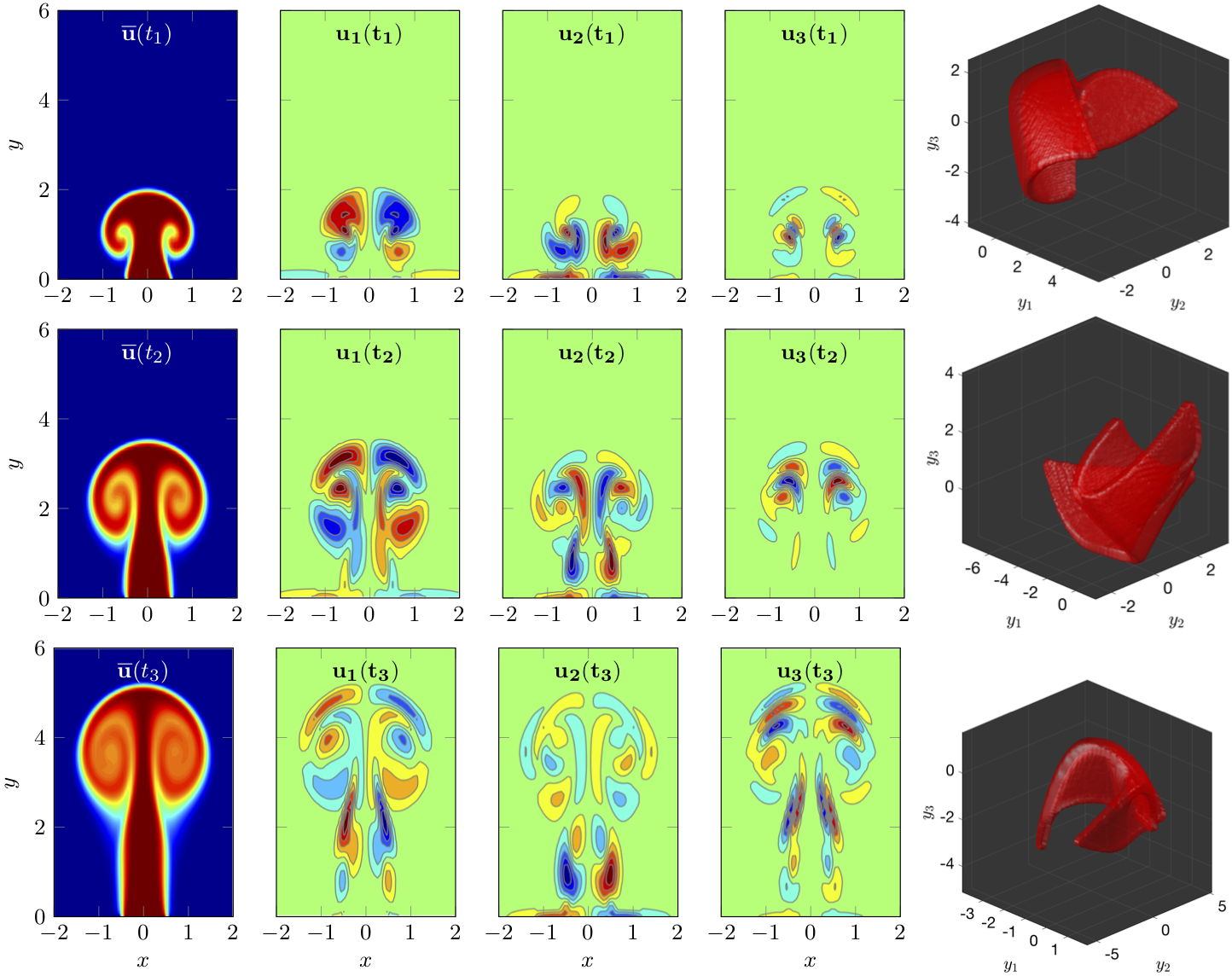}
\caption{Dynamic basis reduction of transient vertical jet in quiescent flow.  Left: Flow mean visualized by passive scalar. Middle: Three most dominant dynamics basis visualized by vorticity $\omega_z$. Right:   Iso-surface of joint probability density function $\pi_{y_1y_2y_3}$ of the first three dynamic modes. Each row shows the state of the reduction at a fixed time with first row: $t=2$; second row: $t=6$; and third row $t=10$.   }
\label{Fig:jet}
\end{figure}

As a reference reduction, we first computed the proper orthogonal decomposition (POD) modes. Since POD modes are often used in deterministic flows, we explain how these modes are computed in a stochastic setting. To do so,  we first formed the matrix $\bm{A}$ that contains all  mean-subtracted snapshots of all 36 samples, i.e. the entire DNS data.  We then performed  singular    value decomposition of this matrix. The POD modes are  then computed as the left singular eigenvectors of matrix $\bm{A}$. Therefore, the computed POD modes are an optimal basis for the entire DNS data in the second norm sense.  In the middle three panels of Figure \ref{Fig:Jet_POD}, the first three POD modes associated with largest singular values are shown.

We solved the  evolution Eqs. \ref{eq:u_ev_DO} and \ref{eq:y_ev_DO} with $r=4$ with fourth order Runge-Kutta integration scheme.  The time snapshots of matrix $\bm{T}$ were obtained by skipping  every three time steps of the DNS, i.e. the  time interval of dynamic basis evolution is $\Delta t= 2\times 10^{-2}$ .   

In Figure \ref{Fig:jet}, we show the results of the state of the dynamics basis for three time snapshots at $t_1=2$, $t_2=6$ and $t_3=10$ time units, with each row showing the results at a fixed time. For the purpose of visualizing the mean, we solved a passive scalar transport equation along with Navier-Stokes equation with Schmidt number $Sc=1$. See \cite{Babaee_PRSA} for more details for similar approach for visualization.  In the leftmost column, the mean data is visualized by the scalar transport intensity.  The mean shows the puff growing in size as the jet  entrains the quiescent background flow while advecting upward. In the middle three columns, the first three most dominant dynamic bases are visualized by their vorticity: $\omega_{z_i} = \partial u_i/\partial y-\partial v_i/\partial x$.  The evolution of the dynamic basis along with the data is clearly seen. The dynamic basis reveals the \emph{spatio-temporal coherent structures} in the flow.  Similar to the POD modes (see Figure \ref{Fig:Jet_POD_Modes}), the higher modes of the dynamic basis  show finer spatial structures.   However, the POD modes are static and they are optimal in a time-average sense, as opposed to the dynamic basis that evolves according to $\bm{\dot{T}}$. We will compare the POD reduction error with the dynamic basis in the section.  In the rightmost column, the low-dimensional attractor is visualized  by plotting the iso-surface of joint probability density function $\pi_{y_1y_2y_3}$ of the three most dominant stochastic coefficients: $\bm{y}_1,\bm{y}_2$ and $\bm{y}_3$. Since the DNS are sampled at PCM points, the polynomial chaos expansion is used to generate more samples for visualizing $\pi_{y_1y_2y_3}$ in the post-processing stage. The attractor  clearly shows a low-dimensional manifold as it evolves in time.  

\begin{figure}[h]
\centering
\subfigure[]{
\includegraphics[width=.47\textwidth]{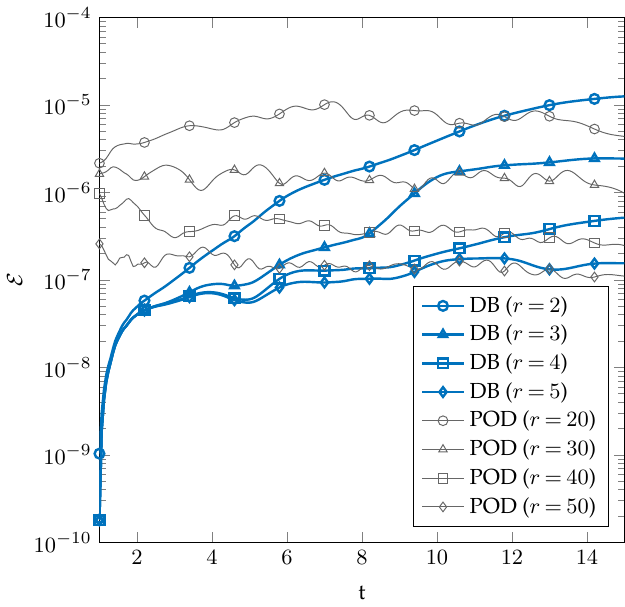}
\label{Fig:Jet_Error}
}
\subfigure[]{
\includegraphics[width=.47\textwidth]{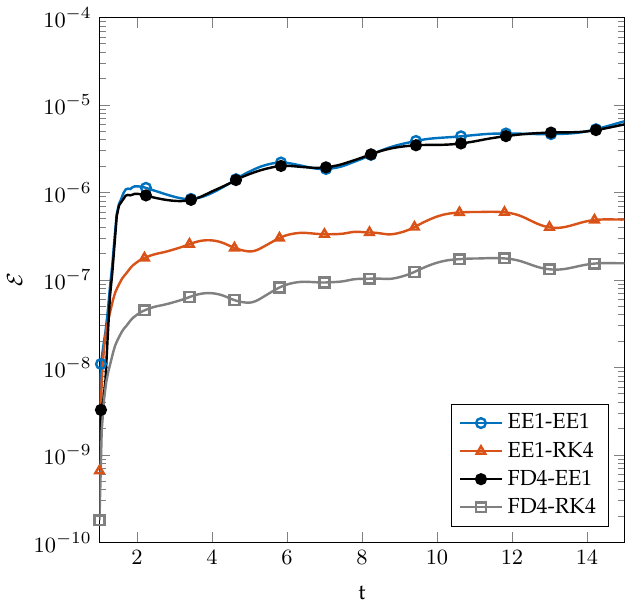}
\label{Fig:Jet_Error_TimeError}
}
\caption{(a) Comparison of the reduction error  between dynamic basis (DB) and proper orthogonal decomposition (POD). (b) Effect of discretization scheme for computing $\bm{\dot{T}}$ and evolution of the dynamic basis. First-order Euler explicit (EE1), fourth-order  central finite difference (FD4) and Runge-Kutta (RK4). Identifiers are shown in the form of X-Y: X shows the discretization scheme used for estimating  $\dot{\bm{T}}$ and Y shows  the time-advancement scheme for solving  evolution equation for the dynamic basis. }
\end{figure}

In Figure \ref{Fig:Jet_Error}, we make a comparison between POD and dynamic basis  reduction errors. This error is obtained by:
\begin{equation*}
\mathcal{E}(t)
 =\left\Vert \mathbb{E}\bigg[ \bm{U} \bm{Y}^T  - \bm{T}  \bigg]^{2}  \right\Vert_F, \label{eq:POD_func}
\end{equation*}
where for POD modes $\bm{U}$ represents the first $r$ dominant POD modes and $\bm{Y}$ is obtained via the orthogonal projection of the data onto $\bm{U}$, i.e.: $\bm{Y}(t)=\inner{\bm{U}}{\bm{T}(t)}$. As shown  in Figure \ref{Fig:Jet_Error}, the reduction of dynamic basis $r=2$ is comparable to that of POD with $r=20$ for long time, and in the early evolution (i.e. $t<4$), the error of dynamic basis ($r=2)$ is smaller than POD with $r=50$. Similar behaviour is observed for dynamic basis $r=3,4$ and 5, where   a ten-fold increase in POD reduction size is required to have a comparable reduction error with dynamic basis. We also observe that the dynamic basis error increases with time. This is due two effects: (1) as the jet advects upward, the uncertainty grows and more modes are required to better estimate the observations; and (2)  the effect of  unresolved modes in the dynamic basis reduction. In the case of drastic reduction, for example $r=2$, discarding the unresolved modes will lead to error increase over time. The effect of unresolved modes can be properly studied in the  Mori-Zwanzig formalism \cite{CHK02}. As it can be seen in Figure \ref{Fig:Jet_Error} increasing the number of dynamic basis significantly decreases the error and the growth rate of the error.

In Figure \ref{Fig:Jet_Error_TimeError}, we show the results of different time discretization schemes for: (1) computing $\dot{\bm{T}}$ with two different schemes of first-order Euler explicit (EE1)  and fourth-order  central finite difference (FD4)  and (2) the time advancement of $\bm{U}(t)$ and $\bm{Y}(t)$ via Equations \ref{eq:u_ev_DO} and \ref{eq:y_ev_DO} with EE1 and fourth-order Runge-Kutta (RK4). All four combinations of the above schemes are considered with the identifiers in the form of X-Y, where the first part of the identifier (X) shows the discretization scheme used for estimating  $\dot{\bm{T}}$ and the second part of the identifier (Y) shows  the time-advancement scheme for solving  Equations \ref{eq:u_ev_DO} and \ref{eq:y_ev_DO}.  In particular we considered these schemes:  EE1-EE1, EE1-RK4, FD4-EE1 and FD4-RK4. We used dynamic basis with the size of $r=5$ for all these schemes. The error of schemes EE1-EE1 and FD4-EE1 are very close -- suggesting that the time-advancement error of solving  solving Equations \ref{eq:u_ev_DO} and \ref{eq:y_ev_DO} is the primary source of error and increasing the accuracy of estimating $\dot{\bm{T}}$ does not improve the overall accuracy. Comparing the errors of EE1-EE1 with EE1-RK4 confirms that the error can be reduced by increasing the accuracy of time integrator of the evolution equations for $\bm{U}$ and $\bm{Y}$. However, comparing the errors of  EE1-RK4 with FD4-RK4 demonstrates that the primary contributor of error for the scheme EE1-RK4 is the estimation of $\dot{\bm{T}}$. The totality of these observations show that  both the accuracy of estimating  $\dot{\bm{T}}$ and solving the evolution equations of the dynamic basis affect the reduction error potentially significantly. Obviously, if  $\Delta t$ is decreased, the  contribution of time-discretization error will become less significant compared to reduction error.

\section{Conclusion}
The objective of this paper is to develop a reduced description  of the transient dynamics that relies on the observables of the dynamics. To this end, we present a variational principle whose optimality conditions lead to closed-form evolution equation for a set of time-dependent orthonormal modes and their associated coefficients, referred to as dynamic basis. Dynamic basis is driven by  the time derivative of the random observations, and it does not require the underlying dynamical system that generates the observables. In that sense, the proposed reduction is purely observation-driven and may apply to systems whose models are not known. 

The application of the dynamic basis was demonstrated for a variety of examples: (1)  advection equation with random wave speed for which exact solution for the dynamic basis was derived and it was used for validation purposes. (2)  Kuramoto-Sivashisky equation subject to random initial condition where it was shown that the dynamic basis  captures the transient intrinsic instability - even in the presence of  high-dimensional random initial condition. The performance of the  dynamic basis reduction was compared against the dynamic mode decomposition. (3) Navier-Stokes equation subject to random boundary condition for which we compared the effectiveness of the dynamic basis versus proper orthogonal decomposition, as an example of static-basis reductions.

The dynamic basis evolution equation has a linear computational complexity with respect to the number of observables and the number of random realizations. This is a favorable feature if applications of large dimension are of interest, for example weather forecast, ocean models and turbulent flows. 


\section*{Acknowledgement}
The authors would like to thank Prof. Karniadakis for useful comments. The author  gratefully acknowledge the allocated computer time on the Stampede supercomputers, awarded by the XSEDE program, allocation number ECS140006. The author has been supported through a NASA Grant number 80NSSC18M0150, and an award by American Chemical Society, Grant number 59081-DN19.

\appendix
\section*{Appendix A. Optimality condition of the variational principle}
In the following we show that the first-order optimality condition of the variational principle \ref{Eqn:min_princ} leads to closed-form evolution equations for $\bm{U}$ and $\bm{Y}$. We first observe that:
\begin{align*}
 \mathcal{G}(\dot{\bm{U}}, \dot{\bm{Y}}) &= \lp \dot{\bm{u}}^T_i \Wx \dot{\bm{u}}_j \rp  \lp \bm{y}^T_i \Wxi \bm{y}_j \rp\\ 
                                         &+ \lp \bm{u}^T_i \Wx \bm{u}_j \rp \lp \dot{\bm{y}}^T_i \Wxi \dot{\bm{y}}_j\rp\\
                                         &+2\lp\dot{\bm{u}}^T_i \Wx \bm{u}_j\rp \lp\bm{y}^T_i \Wxi \dot{\bm{y}}_j\rp\\
                                         &-2 \dot{\bm{u}}^T_i \Wx \dot{\bm{T}} \Wxi \bm{y}_i \\
                                         &-2 \bm{u}^T_i \Wx \dot{\bm{T}} \Wxi \dot{\bm{y}}_i \\
                                         &+ \bigg\|\mathbb{E}[\dot{\bm{T}}]^2\bigg\|_F\\
                                         &+\lambda_{ij} \lp \bm{u}^T_i \Wx \dot{\bm{u}}_j  - \bs{\phi}_{ij} \rp.
\end{align*}
First order optimality condition requires $\pfrac{\mathcal{G}}{\dot{\bm{u}}_k }=0$ and $\pfrac{\mathcal{G}}{\dot{\bm{y}}_k }=0$ for $k=1,2,\dots, r$. We first observe that the gradient of the $\mathcal{G}$ with respect to  $\dot{\bm{u}}_k$ can be expressed as:
\begin{align}\label{eq:sen_u}
\pfrac{\mathcal{G}}{\dot{\bm{u}}_k } &= 2 \Wx \dot{\bm{u}}_j   \lp \bm{y}^T_k \Wxi \bm{y}_j \rp\\ \nonumber 
                                     &+ 2 \Wx      \bm{u}_j    \lp \bm{y}^T_k \Wxi \dot{\bm{y}}_j \rp\\ \nonumber 
                                     &- 2 \Wx \dot{\bm{T}} \Wxi    \bm{y}_k\\ \nonumber                                           &+ \lambda_{jk}\Wx \bm{u}_j \\ \nonumber 
                                     & = \bm{0}.
\end{align}
To eliminate $\lambda_{kj}$, we take the inner product of  the above equation  by $\bm{u}_l$:
\begin{align*}
\bm{u}_l^T \pfrac{\mathcal{G}}{\dot{\bm{u}}_k } &= 2 \bs{\phi}_{lj}  \lp \bm{y}^T_k \Wxi \bm{y}_j \rp\\
                                     &+ 2 \delta_{lj}    \lp \bm{y}^T_k \Wxi \dot{\bm{y}}_j \rp\\
                                     &- 2 \bm{u}_l^T \Wx \dot{\bm{T}} \Wxi    \bm{y}_k\\                                            &+ \lambda_{jk} \delta_{lj}\\                                     &=0.
\end{align*}
where we have made use of the definition of $\bs{\phi}_{lj}(t): = \inner{\bm{u}_l}{\dot{\bm{u}}_j} $. Rearranging the above equation results in:
\begin{equation*}
 \lambda_{lk} = -2 \bs{\phi}_{lj}  \lp \bm{y}^T_k \Wxi \bm{y}_j \rp
   			    - 2 \lp \bm{y}^T_k \Wxi \dot{\bm{y}}_l \rp
 			    + 2 \bm{u}_l^T \Wx \dot{\bm{T}} \Wxi    \bm{y}_k.
\end{equation*}
After multiplying  equation \ref{eq:sen_u} from left by  $\dfrac{\Wx^{-1}}{2}$ and substituting $ \lambda_{jk}$ from the above equation we obtain:
\begin{equation*}
  \dot{\bm{u}}_j   \lp \bm{y}^T_k \Wxi \bm{y}_j \rp = 
  \lp \bm{I} -   \bm{u}_j\bm{u}_j^T \Wx \rp  \dot{\bm{T}} \Wxi \bm{y}_k +
 \bm{u}_j \bs{\phi}_{jl} \lp \bm{y}^T_k \Wxi \bm{y}_l \rp.
\end{equation*}
In the vectorized form the above equation can be written as:
\begin{equation*}
  \dot{\bm{U}} = 
  \lp \bm{I} -   \bm{U}\bm{U}^T \Wx \rp  \dot{\bm{T}} \Wxi \bm{Y} \bm{C}^{-1} +
 \bm{U} \bs{\phi}. 
\end{equation*}
where $\bm{C}_{ij}=\bm{y}_i^T\Wxi \bm{y}_j$ is the covariance matrix.  Similarly, 
\begin{align}\label{eq:sen_y}
\pfrac{\mathcal{G}}{\dot{\bm{y}}_k } &= 2 \Wxi \dot{\bm{y}}_j   \lp \bm{u}^T_k \Wx \bm{u}_j \rp\\ \nonumber 
                                     &+ 2 \Wxi      \bm{y}_j    \lp \dot{\bm{u}}^T_j  \Wx \bm{u}_k \rp\\ \nonumber 
                                     &- 2 \Wxi \dot{\bm{T}}^T \Wx    \bm{u}_k\\ \nonumber           
                                     & = \bm{0}.
\end{align}
In the above expression we have made use of  $\lp\dot{\bm{u}}^T_i \Wx \bm{u}_j \rp \lp\bm{y}^T_i \Wxi \dot{\bm{y}}_j\rp = \lp \bm{y}^T_i \Wxi \dot{\bm{y}}_j\rp \lp\dot{\bm{u}}^T_i \Wx \bm{u}_j\rp $ since both   $\lp\dot{\bm{u}}^T_i \Wx \bm{u}_j \rp$ and $\lp \bm{y}^T_i \Wxi \dot{\bm{y}}_j\rp$ are scalars.
By multiplying the above equation \ref{eq:sen_y} by $\dfrac{\Wxi^{-1}}{2}$ from left and replacing $\bm{u}^T_k \Wx \bm{u}_j = \delta_{kj}$ and $\dot{\bm{u}}^T_j  \Wx \bm{u}_k = \bs{\phi}_{kj}=-\bs{\phi}_{jk}$, we obtain:
\begin{equation*}
  \dot{\bm{y}}_k = 
                  \dot{\bm{T}}^T \Wx    \bm{u}_k 
                  +\bm{y}_j\bs{\phi}_{jk}.
\end{equation*}
In the vectorized form the above equation can be expressed as:
\begin{equation*}
  \dot{\bm{Y}} = 
                  \dot{\bm{T}}^T \Wx    \bm{U}
                  +\bm{Y}\bs{\phi}.
\end{equation*}

\section*{Appendix B. Proof of equivalence of reductions }
\begin{theorem}
 Suppose that  $\{\bm{Y} , \bm{U}  \}$ and   $\{\widetilde{\bm{Y}} , \widetilde {\bm{U}}  \}$
 satisfy the evolution equations (\ref{eq:u_ev} and \ref{eq:y_ev}) with different
choices of $\bs{\phi}$ functions denoted by $\bs{\phi} \in \mathbb{R}^{r
\times r}$ and $\widetilde{\bs{\phi}} \in \mathbb{R}^{r
\times r}$, respectively.
We also assume that the two reductions are initially equivalent, \emph{i.e.
} $\bm{U}(0)=\widetilde {\bm{U}}(0)\bm{R}_0$ and $\bm{Y}(0)=\widetilde {\bm{Y}}(0)\bm{R}_0$, where $R_0 \in \mathbb{R}^{r \times r}$ is an orthogonal rotation matrix. Then $\{\bm{Y} , \bm{U}  \}$ and  $\{\widetilde{\bm{Y}} , \widetilde {\bm{U}}  \}$
 are equivalent for $t>0,$ with a rotation matrix $\bm{R}(t) \in \mathbb{R}^{r
\times r}$ governed by the  matrix differential
equation $\dot{\bm{R}}=\bm{R} \bs{\bs{\phi}}-\widetilde{\bs{\phi}} \bm{R}$.
\end{theorem}
\textbf{proof:} We replace $\bm{U}=\widetilde{\bm{U}}\bm{R}$  and $\bm{Y}=\widetilde{\bm{Y}}\bm{R}$ into the evolution equations (\ref{eq:u_ev} and \ref{eq:y_ev}). The evolution of the dynamic basis can be written as:
\begin{align*}
\dot{\bm{U}} &= \dot{\widetilde{\bm{U}}}\bm{R} + \widetilde{\bm{U}}\dot{\bm{R}}\\
             & = \lp \bm{I} -   \widetilde{\bm{U}}\widetilde{\bm{U}}^T \Wx \rp  \dot{\bm{T}} \Wxi \widetilde{\bm{Y}} \bm{R} \bm{C}^{-1} +
 \widetilde{\bm{U}} \bm{R} \bs{\phi}.
\end{align*}
Therefore: 
\begin{align*}
\dot{\widetilde{\bm{U}}} &= \lp \bm{I} -   \widetilde{\bm{U}}\widetilde{\bm{U}}^T \Wx \rp  \dot{\bm{T}} \Wxi \widetilde{\bm{Y}} \bm{R} \bm{C}^{-1} \bm{R}^T + \widetilde{\bm{U}} ( \bm{R}\bs{\phi} \bm{R}^T -   \dot{\bm{R}} \bm{R}^T ) \\
                        & = \lp \bm{I} -   \widetilde{\bm{U}}\widetilde{\bm{U}}^T \Wx \rp  \dot{\bm{T}} \Wxi \widetilde{\bm{Y}}  \widetilde{\bm{C}}^{-1} + \widetilde{\bm{U}} \widetilde{\bs{\phi}}
\end{align*}
where $\widetilde{\bm{C}} = \mathbb{E}[\widetilde{\bm{Y}}^T \widetilde{\bm{Y}}] = \bm{R} \bm{C} \bm{R}^T  $ and  therefore, $\widetilde{\bm{C}}^{-1} = \bm{R} \bm{C}^{-1} \bm{R}^T$, where $\bm{C}$ are $\widetilde{\bm{C}}$ are similar matrices and have the same eigenvalues. In the last equation we have used $\dot{\bm{R}}=\bm{R} \bs{\bs{\phi}}-\widetilde{\bs{\phi}} \bm{R}$, and made use of the  identity $\bm{R}^{-1}=\bm{R}^T$.
The evolution of the stochastic coefficients can be written as:
\begin{align*}
\dot{\bm{Y}} &= \dot{\widetilde{\bm{Y}}}\bm{R} + \widetilde{\bm{Y}}\dot{\bm{R}}\\
             & =\dot{\bm{T}}^T \Wx    \dot{\widetilde{\bm{U}}}\bm{R}
                  +\dot{\widetilde{\bm{Y}}}\bm{R}\bs{\phi}.
\end{align*}
Therefore: 
\begin{align*}
\dot{\widetilde{\bm{Y}}} &=\dot{\bm{T}}^T \Wx    \dot{\widetilde{\bm{U}}}
                  +\dot{\widetilde{\bm{Y}}}( \bm{R}\bs{\phi}\bm{R}^T- \dot{\bm{R}} \bm{R}^T)\\
                         &=\dot{\bm{T}}^T \Wx    \dot{\widetilde{\bm{U}}}
                       +\dot{\widetilde{\bm{Y}}}\widetilde{\bs{\phi}}.
\end{align*}
This completes the proof. 

\bibliographystyle{unsrt}
\pagebreak
\bibliography{biblio}

\end{document}